\documentclass[12pt,a4paper]{amsart}
\usepackage{amscd}
\usepackage{amsmath}
\usepackage{latexsym}
\usepackage{amsfonts}
\usepackage{amssymb}
\usepackage{amsthm}
\usepackage{graphicx}
\usepackage{verbatim}
\usepackage{mathrsfs}
\usepackage{enumerate}
\usepackage{color}

\newtheorem{definition}{Definition}[section]
\newtheorem{theorem}{Theorem}[section]
\newtheorem{lemma}[theorem]{Lemma}

\newtheorem{proposition}[theorem]{Proposition}

\theoremstyle{remark}

\newcommand{\ZZ}{\mathbb{Z}}

\newcommand{\PP}{\mathbb{P}}

\newcommand{\RR}{\mathbb{R}}

\DeclareMathOperator{\Div}{div}

\allowdisplaybreaks
\begin{document}
\title[Local wellposedness for the non-resistive MHD Equations ]{Local wellposedness for the non-resistive MHD Equations in optimal Sobolev spaces}

\author{Yatao Li}
\address[Y. Li]{The Graduate School of China Academy of Engineering Physics,P.O.Box 2101,Beijing 100088,PR China}
\email{liytao\_maths@163.com}
\date{\today}

\begin{abstract}
 In this paper, we consider the Cauchy problem of the non-resistive magnetohydrodynamics equations in $\RR^d$ for $d=2,3$. We show that the system is  locally well-posed in $H^{s-1}\times H^s$  by establishing a new commutator estimate and utilizing the heat smooth effect in Chemin-Lerner frame.  The space $H^{s-1}\times H^s$ is optimal in Sobolev spaces for the local well-posedness of the system in the scaling sense. Therefore, we improve the results in  \cite{Fef2017}.
\end{abstract}

\maketitle
\section{Introduction}\label{INTR}
\setcounter{section}{1}\setcounter{equation}{0}

In this paper, we consider the Cauchy problem of the following incompressible non-resistive magnetohydrodynamics equations (NMHD) for $d=2,3$:
\begin{equation}\label{eq.NMHD}
\left\{\begin{array}{ll}
\partial_t u-\nu\Delta u+u\cdot\nabla u+\nabla \pi=b\cdot\nabla b,\\
\partial_t b+u\cdot\nabla b=b\cdot\nabla u,\\
\Div u=\Div b=0,
\end{array}\right.
\end{equation}
where vector fields $u=(u^1,u^2,\dots,u^d)$, $b=(b^1,b^2,\dots,b^d)$ are the fluid velocity and the magnetic field respectively, the scalar function $\pi$ is the fluid pressure, and $\nu>0$ is the viscosity coefficient. System (1.1) describes the dynamics of magnetic field in electrically conducting fluid. It has been extensively investigated by mathematicians in the last few decades. We can refer to \cite{CHMDRR16,ChMZh07,CHMZH08,DJL72,JN06,LinXUZH15,STemamR83}.\\
\indent In this paper, we concern the problem of local well-posedness in optimal Sobole spaces. Fefferman et al. obtained local-in-time existence of strong solutions to \eqref{eq.NMHD} in $\RR^d,d=2,3$ with $(u_0,b_0)\in H^{s}\times H^s$ in \cite{Fefet2014} and $(u_0,b_0)\in H^{s-1+\varepsilon}\times H^s$ in \cite{Fef2017}. The aim of this paper is to remove $\varepsilon$ in \cite{Fef2017} and thus obtain local well-posedness in the optimal Sobolev space based on the natural scaling of system (1.1).  The main difficulty comes from nonlinear terms in the transport equation due to the lack of the diffusion of space variable. However,  by applying the frequency localization method and some harmonic analysis techniques, establishing a new commutator estimate and utilizing sufficiently the heat smooth effect in Chemin-Lerner Besov spaces, we overcome the disadvantage. Now, let's state our main result as follows:
\begin{theorem}\label{THM.INT}
 Assume that the initial data that $u_0 \in H^{s-1},\,b_0 \in H^{s},\,s>d/2, \,d=2,3.$ Then there exists a strictly positive maximum time $T_{\ast}$ such that a unique solution $(u,b)$ of the system \eqref{eq.NMHD} exists in the space $C([0,T_{\ast});H^{s-1}\times H^{s}).$ Moreover, the solution
$u\in  \widetilde L^2([0,T_{\ast});B^{s}_{2,2})\cap \widetilde L^1([0,T_{\ast});B^{s+1}_{2,2}).$
\end{theorem}
\indent The paper is organized as follows. In Section $2$, we recall Littlewood-Paley theory and give some properties of Besov space. In Section $3$, we prove the local existence and uniqueness of the solution of system (\ref{eq.NMHD}).\\
$\mathbf{Notation}.$ We denote by $\langle\cdot,\cdot\rangle$ the inner product on $\RR^d$. Given a Banach space $X$, we denotes its norm by $\|\cdot\|_X$. The uniform constant $C$ may be different on different lines in this paper.
\section{Preliminaries}\label{PRE}
\setcounter{section}{2}\setcounter{equation}{0}
In this section, we firstly recall some Littlewood-Paley theory.  One can refer to \cite{BCD11,MWZ12} for more details.\\
\indent Let $\varphi,\,\chi\in \mathcal{S}(\RR^d)$ be two smooth radial functions with values in $[0,1]$. $\varphi$ is supported in the annulus $\{\xi\in \RR^d:\frac{3}{4}\le |\xi|\le\frac{8}{3}\}$, $\chi$ is supported in the  unit ball $\mathcal B(0,1)$ of $\RR^d$. They satisfy with
\begin{center}
$\chi(\xi)+\sum_{j\ge0}\varphi_j(\xi)=1$\quad\quad\quad for\quad$\xi\in\RR^d,$
\end{center}
\begin{center}
$\sum_{j\in\ZZ}\varphi_j(\xi)=1$\quad\quad\quad for\quad $\xi\in\RR^d \setminus \{0\},$
\end{center}
where we denote $\varphi_j(\xi)=\varphi(2^{-j}\xi)$.\\
Let us denote the Fourier transform on $\RR^d$ by $\mathcal{F}$ and write $h=\mathcal{F}^{-1}{\varphi}$ and $\widetilde {h}=\mathcal{F}^{-1}{\chi}$.  The homogeneous localization operator $\dot{\Delta}_j$ and the homogeneous low-frequency cut-off operators $ \dot{S}_j$ are defined for all $j\in\ZZ$ by
$$\dot{\Delta}_ju=\varphi_j(D)u=2^{jd}\int_{\RR^d}h(2^jy)u(x-y)\,dy, \;\;$$
$$ \dot{S}_ju=\sum_{k\le j-1}\dot{\Delta}_ku=2^{jd}\int_{\RR^d}\widetilde h(2^jy)u(x-y)\,dy,$$
and the inhomogeneous localization oprator:
\begin{center}
$\Delta_ju=\varphi_j(D)u=\mathcal{F}^{-1}(\varphi_j(\xi)\hat{u}),$\quad \;\,if \; \quad $j\ge0;$
\end{center}
\begin{center}
 $\Delta_{-1}u=\mathcal{F}^{-1}(\chi(\xi)\hat{u});$\quad $\Delta_ju=0$,\quad if \quad $j\le-2.$
\end{center}
The inhomogeneous low-frequency cut-off operators $ {S}_j$ are defined for all $j\in\ZZ$ by
\begin{center}
$ S_ju=\sum_{k\le j-1}\Delta_ku=\mathcal{F}^{-1}(\chi(2^{-j}\xi)\hat{u}).$
\end{center}
\noindent One can easily verify that
\begin{center}
$\dot{\Delta}_j\dot{\Delta}_{j'}u=0$ \quad\quad\quad\quad  if \quad $|j-j'|\ge2,$\\
$\dot{\Delta}_j(\dot{S}_{j'-1}u\dot{\Delta}_{j'}u)=0$\quad if \quad$|j-j'|\ge5.$
\end{center}

\indent Next, we recall Bony's decomposition from \cite{BCD11}:
\begin{definition} \label{defparaproduct}
For any  $u, v\in\mathcal {S}^{\prime}/{\mathcal P}(\RR^d)$, $u
v$ has the homogeneous Bony paraproduct decomposition:
$$u v=\dot{T}_{u} v+\dot{T}_{v}u+\dot{R}(u,v),
$$
where $$\dot T_{u}v:=\sum_{j\leq{k-2}}\dot\Delta_{j}u\dot\Delta_{k}v=\sum_{j}{\dot{S}_{j-1}u}{\dot{\Delta}_{j}v}  \quad and \quad \dot{R}(u,v)=:\sum_{|j-j'|\leq1}{\dot\Delta}_{j}u{\dot\Delta}_{j'}v.$$
For any  $u, v\in\mathcal {S}^{\prime}$,
$uv$ has the inhomogeneous Bony paraproduct decomposition:
$$u v={T}_{u} v+{T}_{v}u+{R}(u,v),
$$
where $$ T_{u}v:=\sum_{j\leq{k-2}}\Delta_{j}u\Delta_{k}v=\sum_{j}{{S}_{j-1}u}{{\Delta}_{j}v}  \quad and \quad {R}(u,v)=:\sum_{|j-j'|\leq1}{\Delta}_{j}u{\Delta}_{j'}v.$$
\end{definition}
\indent We will use repeatedly the following classical Bernstein-Type lemma
\begin{lemma}\cite{BCD11}\label{lem.Ber}
Let $\mathcal{C}$ be an annulus, $\mathcal B$ a ball, and $(p,q)\in [1,\infty]^2$ with $1\leq p\leq q$. Then for any vector filed $f\in L^p(\RR^d)$, there exist a constant $C>0$, independent of $f$ and $\lambda$, such that for any $k\in \ZZ$,
\[\|D^k f\|_{L^q}\leq C^{k+1}\lambda^{k+d(\frac1p-\frac1q)}\|f\|_{L^p}\quad \text{if}\quad {\mathop{\rm supp}\nolimits}\, \hat{f}\subset \lambda \mathcal B,\]
\[C^{-k-1}\lambda^{k}\|f\|_{L^p}\leq \|D^k u\|_{L^p}\leq C^{k+1}\lambda^{k}\|f\|_{L^p}\quad \text{if}\quad {\mathop{\rm supp}\nolimits}\, \hat{f}\subset \lambda \mathcal{C}.\]
\end{lemma}
The definition of Besov spaces is as follows:
\begin{definition}\label{Besovdef}
Let $s\in \RR$, $1\le p,r \le \infty$. $S'$ be the space of tempered distributions and $\mathcal{P}$ is the set of all polynomials. The homogeneous  Besov space are defined as follows
$$\dot{B}^{s}_{p,r}(\RR^d):=\{u\in \mathcal{S}'(\RR^d)/\mathcal{P}:\|u\|_{\dot{B}^{s}_{p,r}(\RR^d)}<\infty\},$$
where
$$\|u\|_{\dot{B}^{s}_{p,r}(\RR^d)}=(\sum\limits_{j\in \ZZ}2^{js}\|\dot{\Delta}_ju\|^{r}_{L^{p}(\RR^d)})^{\frac{1}{r}};\quad\quad\quad\quad$$
The inhomogeneous  Besov space are defined as follows
$${B}^{s}_{p,r}(\RR^d):=\{u\in \mathcal{S}'(\RR^d):\|u\|_{{B}^{s}_{p,r}(\RR^d)}<\infty\},\quad$$
where
$$\|u\|_{{B}^{s}_{p,r}(\RR^d)}=(\sum\limits_{j\ge {-1}}2^{js}\|{\Delta}_ju\|^{r}_{L^{p}(\RR^d)})^{\frac{1}{r}}.\quad\quad\quad$$
\end{definition}
\indent $Remarks$.\quad  When $p=r=2$, let us point out that for any $s\in \RR$, ${\dot B}^s_{2,2}$ and ${ B}^s_{2,2}$ are the usual Sobolev space ${\dot H}^s$ and ${H}^s$, respectively. In addition, ${ B}^{\delta}_{2,2}={ \dot B}^{\delta}_{2,2}\cap L^2$ with $\delta\ge0$. \\
\indent We also need use the Chemin-Lerner type homogeneous Besov space(see \cite{BCD11,MWZ12}):
\begin{definition}\label{Besovt}
Let $s\in \RR$, $1\le p,q,r \le \infty$, and $T\in(0,\infty]$. the time homogeneous Besov space $\widetilde {L}^q_T{\dot B}^{s}_{p,r}(\RR^d)$ are defined as follows
\begin{equation*}
\|u\|_{\widetilde {L}^q_T({\dot B}^{s}_{p,r}(\RR^d))}:=(\sum\limits_{j\in \ZZ}2^{jsr}\|\Delta_ju\|^{r}_{ {L}^q_T(L^{p}(\RR^d))})^{\frac{1}{r}}.
\end{equation*}
\end{definition}
\indent By Minkowski's inequality, it is easy to get that:
$$\|u\|_{\widetilde {L}^q_T({\dot B}^{s}_{p,r}(\RR^d))}\leq \|u\|_{ {L}^q_T({\dot B}^{s}_{p,r}(\RR^d))}\quad if \quad q\leq r,$$ $$\|u\|_{\widetilde {L}^q_T({\dot B}^{s}_{p,r}(\RR^d))}\ge \|u\|_{{L}^q_T({\dot B}^{s}_{p,r}(\RR^d))}\quad if \quad q\ge r.$$
The inhomogeneous case is similar (see \cite{BCD11,MWZ12}).\\
\indent Some useful properties of the Besov spaces or the Chemin-Lerner type Besov space from \cite{Danchin01} are collected as follows:
\begin{lemma}\label{Besovpro}
For all $s,s_1,s_2\in\mathbb{R},1\leq p,r,r_1,r_2,q_1,q_2\leq+\infty,\frac{1}{r}\le\frac{1}{r_1}+\frac{1}{r_2}\le1\,$ and $ \frac{1}{q}=\frac{1}{q_1}+\frac{1}{q_2}$,   \\
(i)if $s_1,s_2\leq\frac{d}{p}$ such that $s_1+s_2> d \max\{0,\frac 2p-1\}$, $u\in{\dot B}^{s_1}_{p,r_1}$ and $v\in {\dot B}^{s_2}_{p,r_2}$. Then there hold that
$$\|uv\|_{{\dot B}^{s_1+s_2-\frac{d}{p}}_{p,r}}\le C\|u\|_{{\dot B}^{s_1}
_{p,r_1}}\|v\|_{{\dot B}^{s_2}_{p,r_2}};\quad\quad\quad\quad\quad$$
$$\|uv\|_{\widetilde {L}^q_T({\dot B}^{s_1+s_2-\frac{d}{p}}_{p,r})}\le C\|u\|_{\widetilde {L}^{q_1}_T({\dot B}^{s_1}
_{p,r_1})}\|v\|_{\widetilde {L}^{q_2}_T({\dot B}^{s_2}_{p,r_2})};$$
(ii)if $s>0$,\;$\|uv\|_{{\dot B}^{s}_{p,r}}\leq C\|u\|_{L^{\infty}}\|v\|_{{\dot B}^{s}_{p,r}}+\|u\|_{{\dot B}^{s}
_{p,r}}\|v\|_{L^{\infty}};$\\
(iii) if $p_1\leq p_2$, $r_1\leq r_2$, then $$\dot{B}^s_{p_1,r_1}\hookrightarrow \dot{B}^{s-\frac{d}{p_1}+\frac{d}{p_2}}_{p_2,r_2},\quad{B}^s_{p_1,r_1}\hookrightarrow{B}^{s-\frac{d}{p_1}+\frac{d}{p_2}}_{p_2,r_2};\quad\quad$$
(iv) if $s_1\neq s_2$ and $\theta\in(0,1)$, then
$$\|u\|_{{B}_{p,r}^{\theta s_1+(1-\theta)s_2}}\leq\|u\|^\theta_{{B}^{s_1}_{p,r}}\|u\|^{1-\theta}_
{{B}^{s_2}_{p,r}}.\quad\quad\quad\quad\quad\quad$$
\end{lemma}
\indent We will present some estimates for the heat equation
\begin{align}\label{eq.heat}
\partial_t f-\nu\Delta f=g,\;f|_{t=0}=f_0,
\end{align}
in homogenous Besov spaces (see \cite{BCD11,MWZ12}).
\begin{lemma}\label{pro.heat}
Let $\rho,\,\rho_1,\,p,$\;and\;$r$\;satisfy that\;$1\leq p,\,r\leq{\infty}$,\;and\;$1\leq \rho_1\leq\rho\leq\infty.$\; $f_0\in {\dot B}^{s}_{p,r},\; g\in {{\widetilde L}^{\rho_1}_T({\dot B}^{s-2+\frac 2{\rho_1}}_{p,r})}$ and $f$ is a solution of equation \eqref{eq.heat}. Then there exists an absolute constant $C$ such that
$${\nu}^{\frac 1\rho}\|f\|_{{\widetilde L}^{\rho}_T({\dot B}^{s+\frac 2\rho}_{p,r})}\le C\Big(\|f_0\|_{{\dot B}^{s}_{p,r}}+{\nu }^{\frac 1\rho_1-1}\|g\|_{{\widetilde L}^{\rho_1}_T({\dot B}^{s-2+\frac 2{\rho_1}}_{p,r})}\Big). $$
\end{lemma}
\indent We also use the notation of the commutator:
$$[\dot\Delta_j,v\cdot\nabla]u=\dot\Delta_j(v\cdot\nabla)-v\cdot\nabla \dot\Delta_ju,\quad\quad[\Delta_j,v\cdot\nabla]u=\Delta_j(v\cdot\nabla)-v\cdot\nabla \Delta_ju.$$There are two commutator estimates to be used. One  will be applied in estimates for the transport equation and is as following:
\begin{lemma}\cite{BCD11}\label{Dancom01}
Let $s\in\RR,\,1\leq r\leq\infty,$  $1\leq p\leq p_1\leq\infty,$ and $s<\frac dp$ or $s=\frac dp$ if $r=1$. Let $v$ be a divergence-free vector field over $\RR^d$. Assume that $ -1-d\{\frac 1p_1,\frac 1{p'}\}<s<1+\frac dp_1$,
then $$\|[\Delta_j,v\cdot\nabla]u\|_{ B^s_{p,r}}\leq C\|\nabla v\|_{ B^{\frac dp_1}_{p_1,r}\cap L^{\infty}}\|u\|_{ B^s_{p,r}}.$$
\end{lemma}
\noindent The other is  established by us for the first time, which plays an important role in our proof. We give its detailed proof as follow.
\begin{proposition}\label{pro.comE}
Let $s>0,\;1\leq p,q\leq \infty$, we have the following commutator estimate,
\[\|[\dot{\Delta}_j,\, u\cdot\nabla] v\|_{L^p}\le C c_j2^{-j(s-1)}\big(\|u\|_{L^{\infty}}\|v\|_{\dot{B}^s_{p,q}}+\|v\|_{L^{\infty}}\|u\|_{\dot{B}^s_{p,q}}\big),\]
where $\{c_j\}_{j\in\ZZ}\in {\ell^q}$.
\end{proposition}
\begin{proof}
Using Bony's paraproduct decomposition by Definition \ref{defparaproduct}, we have
\begin{align*}
[\dot{\Delta}_j,\, u\cdot\nabla] v=&[\dot{\Delta}_j,\, \dot{T}_{u^i}]\partial_i v+\dot{\Delta}_j\big(\dot{T}_{\partial_i v} u^i\big)+\dot{\Delta}_j\big(\dot{R}(u^i,\,\partial_i v)\big)-\dot{T}_{\partial_i \dot{\Delta}_jv} u^i-\dot{R}\big(u^i,\,\partial_i \dot{\Delta}_j v\big)
\end{align*}
For  the  term $[\dot{\Delta}_j,\, \dot{T}_{u^i}]\partial_i v$. By the definition of $\dot{\Delta}_j$ and Taylor's formula, we have
\begin{align*}
[\dot{\Delta}_j,\, \dot{T}_{u^i}]\partial_i v=&\sum_{|j-j'|\leq 4}\big(\dot{\Delta}_j(\dot{S}_{j'-1}u^i\dot{\Delta}_{j'}\partial_iv)-
\dot{S}_{j'-1}u^i\dot{\Delta}_{j'}\dot{\Delta}_{j}\partial_iv\big)\\
=&\sum_{|j-j'|\leq 4} 2^{jd}\int_{\RR^d}\varphi(2^jy)\Big(\dot{S}_{j'-1}u^i(x-y)-\dot{S}_{j'-1}u^i(x)\Big)
\dot{\Delta}_{j'}\partial_iv(x-y)\,\mathrm{d}y\\
=&\sum_{|j-j'|\leq 4} 2^{jd}\int_{\RR^d}\int^1_0(-y)\cdot\nabla\dot{S}_{j'-1}u^i(x-\tau y)\,\mathrm{d}\tau\,
\varphi(2^jy)\dot{\Delta}_{j'}\partial_iv(x-y)\,\mathrm{d}y.
\end{align*}
By Minkowski's inequality and Lemma \ref{lem.Ber}, we get that
\begin{equation*}
\|[\dot{\Delta}_j,\, \dot{T}_{u^i}]\partial_i v\|_{L^p}\leq C\sum_{|j-j'|\leq 4}\| \dot{S}_{j'-1} \nabla u^i\|_{L^{\infty}}\|\dot{\Delta}_{j'} v\|_{L^p}\leq C\sum_{|j-j'|\leq 4}2^{j'}\| u\|_{L^{\infty}}\|\dot{\Delta}_{j'} v\|_{L^p}
\end{equation*}
For $\dot{\Delta}_j\big(\dot{T}_{\partial_i v} u^i\big)$ and $\dot{T}_{\partial_i \dot{\Delta}_jv} u^i$, we have
\begin{equation*}
\|\dot{\Delta}_j\big(\dot{T}_{\partial_i v} u^i\big)\|_{L^p}\leq C\sum_{|j-j'|\leq 4}\| \dot{S}_{j'-1}\partial_i v\|_{L^{\infty}}\|\dot{\Delta}_{j'}u^i\|_{L^p}\leq C\sum_{|j-j'|\leq 4}2^{j'}\| v\|_{L^{\infty}}\|\dot{\Delta}_{j'}u\|_{L^p}
\end{equation*}
and
\begin{equation*}
\|\dot{T}_{\partial_i \dot{\Delta}_jv} u^i\|_{L^p}\leq C\sum_{j'> j}\|\dot{\Delta}_j\dot{S}_{j'-1}\partial_i v\|_{L^{\infty}}\|\dot{\Delta}_{j'}u^i\|_{L^p}\leq C\sum_{j'> j}2^j\| v\|_{L^{\infty}}\|\dot{\Delta}_{j'}u\|_{L^p}.\quad\quad\quad
\end{equation*}
For the remainder terms $\dot{\Delta}_j\big(\dot{R}(u^i,\,\partial_i v)\big)$ and $\dot{R}\big(u^i,\,\partial_i \dot{\Delta}_j v\big)$, we have that
\begin{equation*}
\|\dot{\Delta}_j\big(\dot{R}(u^i,\,\partial_i v)\big)\|_{L^p}\leq C2^j\sum_{j'\geq j-3}\|\dot{\Delta}_{j'}u^i\|_{L^{p}}\|\tilde{\dot{\Delta}}_{j'}v\|_{L^{\infty}}\leq C2^j\sum_{j'\geq j-3}\|\dot{\Delta}_{j'}u^i\|_{L^{p}}\|v\|_{L^{\infty}}
\end{equation*}
and
\begin{equation*}
\|\dot{R}\big(u^i,\,\partial_i \dot{\Delta}_j v\big)\|_{L^p}\leq C2^j\sum_{|j'-j|\leq 1}\|\dot{\Delta}_{j'}u^i\|_{L^{p}}\|\tilde{\dot{\Delta}}_{j'}v\|_{L^{\infty}}\leq C2^j\sum_{|j'-j|\leq 1}\|\dot{\Delta}_{j'}u^i\|_{L^{p}}\|v\|_{L^{\infty}}.
\end{equation*}
Collecting all the above estimates, we obtain that
\begin{align*}
&2^{j(s-1)}\|[\dot{\Delta}_j,\, u\cdot] v\|_{L^p}\\
\leq& C\sum_{|j-j'|\leq 4}2^{(j-j')(s-1)}2^{j's}\|\dot{\Delta}_{j'}v\|_{L^p}\| u\|_{L^{\infty}}+C\sum_{|j-j'|\leq 4}2^{(j-j')(s-1)}2^{j's}\|\dot{\Delta}_{j'} u\|_{L^p}\|v\|_{L^{\infty}}\\
&+C\sum_{j'> j}2^{(j-j')s}2^{j's}\|\dot{\Delta}_{j'}u\|_{L^p}\| v\|_{L^{\infty}}+C\sum_{j'\geq j-3}2^{(j-j')s}2^{j's}\|\dot{\Delta}_{j'}u\|_{L^{p}}\|v\|_{L^{\infty}}\\
&+C\sum_{|j'-j|\leq 1}2^{(j-j')s}2^{j's}\|\dot{\Delta}_{j'} u\|_{L^{p}}\| v\|_{L^{\infty}}\\
\leq&C\Big(\sum_{|j-j'|\leq 4}2^{(j-j')(s-1)}+\sum_{j'> j}2^{(j-j')s}+\sum_{j'\geq j-3}2^{(j-j')s}+\sum_{|j'-j|\leq 1}2^{(j-j')s}  \Big)\\
&\times c^1_j\|u\|_{\dot{B}^s_{p,q}}\|v\|_{L^{\infty}}+C\sum_{|j-j'|\leq 4}2^{(j-j')(s-1)}c^2_j\| u\|_{L^{\infty}}\|v\|_{\dot{B}^s_{p,q}}
\end{align*}
where $\{c^1_j\}$, $\{c^2_j\}\in \ell^{q}(\ZZ)$. Due to $s>0$, summing in $j\in\ZZ$ yields the desired result.
\end{proof}
\indent Finally, we state a nonlinear Gronwall's inequality.
\begin{lemma}\cite{BNR15}\label{lem.gron}
Assume $x\in W^{1,1}([0,T])\cap C([0,T])$ such that
\[\dot{x}\leq c(t)x^p+e(t),\quad x(0)=x_0\]
with $p>1$, $c,\,e\in L^1([0,T])$. Then for each $t\in[0,T]$, we have
\[x(t)\leq \big(x_0+\int^t_0 e(\tau)\,\mathrm{d}\tau\big)\Big(1-(p-1)\big(x_0+\int^t_0 e(\tau)\,\mathrm{d}\tau\big)^{p-1}\int^t_0c(\tau)\,\mathrm{d}\tau\Big)^{-\frac1{p-1}}.\]
\end{lemma}
\section{Proof of  Theorem 1.1 } \label{LOC}
\setcounter{section}{3}\setcounter{equation}{0}
In this section, we will prove the existence and uniqueness of the solution to the system \eqref{eq.NSFPR}. Firstly, we apply  the very classical Friedrichs method to construct an  approximate system of \eqref{eq.NMHD} in space $\RR^d$.  Next, because the second equation in the whole system is transport, this leads to a loss of one derivative. Therefore, we establish Proposition \ref{pro.comE} and take full advantage of the heat smooth effect to obtain the uniform bound for approximate solutions in Chemin-Lerner type Besov space. Finally, we prove the strong convergence of the sequence and uniqueness in a weaker norm.  The proof of Theorem \ref{THM.INT} need three steps.
\subsection{Construction of approximate solutions to the system\;$\eqref{eq.NMHD}$}\label{CON}
\quad\quad \\ \indent We shall apply the classical Friedrichs method by cut-off in the frequency space.\\
Define the operator $J_n$ by$$J_nu(x):=\mathcal{F}^{-1}(1_{B(0,n)}\hat{u}(\xi)),$$
where $\mathcal{F}$ denotes the Fourier transform in the space variables.\\
Let us  construct the approximate system of \eqref{eq.NMHD} as follows
\begin{equation}\label{eq.NMHDAP}
\left\{\begin{array}{ll}
\partial_tu_n-\nu J_n\Delta u_n=-J_n\PP[(J_nu_n\cdot\nabla)J_nu_n]+J_n\PP[(J_nb_n\cdot\nabla)J_nb_n],\\
\partial_tb_n+J_n[(J_nu_n\cdot\nabla)J_nb_n]=J_n[(J_nb_n\cdot\nabla)J_nu_n]\\
(u_n,b_{n})|_{t=0}=(J_nu_0,J_nb_{0}).\\
\end{array}\right.
\end{equation}
Define
\begin{center}\label{def.V}
  $V_{n}:=\Big\{(u,b)|(u,b)\in L^2\times L^2,$ $ \hat{u}$ and $\hat{b}$  are all supported  in  $B(0,n),\,
\Div u=\Div b=0$\Big\}
\end{center}
endowed the norm with
\begin{center}
 $\|Z\|^2\stackrel{\mathrm{def}}{=}\|u\|^2_{L^2}+\|b\|^2_{L^2},$ for any $Z=(u,b)\in V_{n}.$
\end{center}
\indent The system \eqref{eq.NMHDAP} turns to be an ordinary differential system in the $V_{n}$. Then, the Cauchy-Lipschitz theorem (see Theorem 3.1 in Majda and Bertozzi []) guarantees there exists a unique solution $(u_n,b_n)$ of the system \eqref{eq.NMHDAP} on $[0,T_n)$ for every fixed $n$  where $T_n$ is strictly positive.
Noting that $J_n^2=J_n$, we find that $(J_nu_n,J_{n}b_n)$ is also a solution to system \eqref{eq.NMHDAP}. By the uniqueness of the solution of ODE system \eqref{eq.NMHDAP}, we have $J_nu_n=u_n,\,J_nb_n=b_n$. So $(u_n,b_n)$ is also the solution of the following system
\begin{equation}\label{eq.NMHDAPp}
\left\{\begin{array}{ll}
\partial_tu_n-\nu \Delta u_n=-J_n\PP[(u_n\cdot\nabla)u_n]+J_n\PP[(b_n\cdot\nabla)b_n],\\
\partial_tb_n+J_n[(u_n\cdot\nabla)b_n]=J_n[(b_n\cdot\nabla)u_n]\\
(u_n,b_{n})|_{t=0}=(J_nu_0,J_nb_{0}).\\
\end{array}\right.
\end{equation}
The solution will continue  provided $\|u_n\|_{H^{s-1}}$ and $ \|b_n\|_{H^{s}}$ remain finite.
\subsection{A priori estimate}\label{UNI}
\indent \\ \indent In this section, we will establish the uniformly bound estimate for the smooth approximate solutions of system \eqref{eq.NMHDAPp},  i.e. the following proposition:
\begin{proposition}\label{pro-pri}
 Under the initial condition  of Theorem \ref{INTR}, there exists two positive constants $C_{\ast}=C_{\ast}(\nu,u_0,\,\|b_0\|_{B^{s}_{2,2}})$ and  $T_{\ast}=T_{\ast}(\nu, u_0,\,\|b_0\|_{B^{s}_{2,2}})$ independent of $n$ such that
\begin{align}\label{bound}
\|u_n\|^2_{\widetilde L^{\infty}_{T_*}(B^{s-1}_{2,2})}+\|u_n\|_{\widetilde L^1_{T_*}(B^{s+1}_{2,2})}+\|u_n\|^2_{\widetilde L^2_{T_*}(B^{s}_{2,2})}+\|b_{n}\|^2_{\widetilde L^{\infty}_{T_*}(B^{s}_{2,2})}
\le C_{\ast}.
\end{align}
\end{proposition}
\begin{proof}
\indent For convenience, we omit the  indexes $n$ and $J_n$ in the system \eqref{eq.NMHDAPp}. Therefore,
we only need to do a priori  estimate for the following system:
\begin{eqnarray}\label{eq.NMHDpri}
\begin{cases}
\partial_tu-\nu\Delta u+(u\cdot\nabla)u+\nabla \pi=(b\cdot\nabla)b,\\
\partial_tb+(u\cdot\nabla)b=(b\cdot\nabla)u,\\
\Div u=\Div b=0,\\
u|_{t=0}=u_0,b|_{t=0}=b_0.
\end{cases}
\end{eqnarray}
\indent Firstly, we do some $L^2$ energy estimates. Multiplying the first equation of system \eqref{eq.NMHDpri} by $u$ and multiplying the second equation of system \eqref{eq.NMHDpri} by $b$, using the divergence free condition $\Div u=0$, integrating by parts, we get
\begin{eqnarray*}
\frac{1}{2}\partial_t(\|u\|_{L^2}^2+\|b\|_{L^2}^2)+\nu\|\nabla u\|_{L^2}^2=0.
\end{eqnarray*}
After integration in time  on $[0, T],$ we have that
\begin{equation}\label{L^2}
\begin{split}
&\|u\|^2_{L^2}+\|b\|^2_{L^2}+\nu\int_0^T\|\nabla u\|_{L^2}^2dt=\|u\|^2_{L^2}+\|b\|^2_{L^2}\leq\|u_0\|^2_{B^{s-1}_{2,2}}+\|b_0\|^2_{B^s_{2,2}})=:M_0.
\end{split}
\end{equation}
\indent Next, applying the frequency localization operator $\dot{\Delta}_j$ to the first equation and ${\Delta}_j$ to the second one of the system \eqref{eq.NMHDpri}, we have that
\begin{equation}\label{eq.loc}
\left\{\begin{array}{ll}
\partial_t\dot{\Delta}_j u-\nu\Delta \dot{\Delta}_j u+\dot{\Delta}_j(u\cdot\nabla u)+\nabla\dot{\Delta}_j \pi=\dot{\Delta}_j (b\cdot\nabla b),\\
\partial_t{\Delta}_j b+{\Delta}_j(u\cdot\nabla b)=0,\\
\Div \dot{\Delta}_j u=0,\\
(\dot{\Delta}_j u,{\Delta}_j b)|_{t=0}=(\dot{\Delta}_j u_0,{\Delta}_j b_0).
\end{array}\right.
\end{equation}
Taking the $L^2$ inner product of the first equation of \eqref{eq.loc} with $\dot\Delta_j u$, we get that, by integrations by parts,
\begin{align*}
\frac12\partial_t\|\dot{\Delta}_j u\|^2_{L^2}+\nu\|\nabla \dot{\Delta}_j u\|^2_{L^2}\leq &-\int_{\RR^d}[\dot{\Delta}_j,\,u\cdot \nabla]u\cdot\dot{\Delta}_j u\,dx+\int_{\RR^d}[\dot{\Delta}_j,\,b\cdot \nabla]b\cdot\dot{\Delta}_j u\,dx\\
\leq& \|[\dot{\Delta}_j,\,u\cdot \nabla]u\|_{L^2}\|\dot{\Delta}_j u\|_{L^2}+\|[\dot{\Delta}_j,\,b\cdot \nabla]b\|_{L^2}\|\dot{\Delta}_j u\|_{L^2}.
\end{align*}
The above  inequality holds since $\Div u=0$ which implies that
\[\int_{\RR^d}(u\cdot \nabla)\dot{\Delta}_j u\cdot\dot{\Delta}_j u\,dx=0,\quad\int_{\RR^d}(b\cdot \nabla)\dot{\Delta}_j b\cdot\dot{\Delta}_j u\,dx=0.\]
Integrating the above inequality over $[0,t]$, by Lemma \ref{lem.Ber}, we have that
\begin{align*}
&\|\dot{\Delta}_j u(t)\|^2_{L^2}+c2\nu 2^{2j}\|\dot{\Delta}_j u\|^2_{L^2([0,t]\times\RR^d)}\\
\leq &\|\dot{\Delta}_j u_0\|^2_{L^2}+\int^t_0\|[\dot{\Delta}_j,\,u\cdot \nabla]u\|_{L^2}\|\dot{\Delta}_j u\|_{L^2}+\|[\dot{\Delta}_j,\,b\cdot \nabla]b\|_{L^2}\|\dot{\Delta}_j u\|_{L^2}\,dt^{\prime}.
\end{align*}
Taking $L^{\infty}([0,t])$ of the above inequality on $t$, then using Lemma \ref{Besovpro} and Proposition \ref{pro.comE}, we deduce that
\begin{align*}
&\|\dot{\Delta}_j u\|^2_{L^{\infty}([0,t];\,L^2)}+c2\nu 2^{2j}\|\dot{\Delta}_j u\|^2_{L^2([0,t]\times \RR^d)}\\
\leq &\|\dot{\Delta}_j u_0\|^2_{L^2}+Cc_j2^{-j(s-1)}\int^t_0\|u\|_{L^{\infty}}\|u\|_{\dot{B}^{s}_{2,2}}
\|\dot{\Delta}_ju\|_{
L^2}dt^{\prime}\\
&+Cc_j
2^{-j(s-1)}\int^t_0\|b\|_{L^{\infty}}\|b\|_{\dot{B}^{s}_{2,2}}\|\dot{\Delta}_ju\|_{
L^2}dt^{\prime}.
\end{align*}
Multiplying both sides of the above inequality by $2^{2j(s-1)}$ and then summing over $j$, we get that
\begin{align*}
&\|u\|^2_{\widetilde{L}^{\infty}([0,t];\,\dot{B}^{s-1}_{2,2})}+c2\nu \|u\|^2_{\widetilde{L}^2([0,t];\,\dot{B}^s_{2,2})}\\
\leq &\|u_0\|^2_{\dot{B}^{s-1}_{2,2}}+C\int^t_0\|u\|_{L^{\infty}}\|u\|_{\dot{B}^{s}_{2,2}}\|u\|_{
\dot{B}^{s-1}_{2,2}}dt^{\prime}+C
\int^T_0\|b\|_{L^{\infty}}\|b\|_{\dot{B}^{s}_{2,2}}\|u\|_{\dot{B}^{s-1}_{2,2}}dt^{\prime}.
\end{align*}
This, together with \eqref{L^2}, yields that
\begin{align*}
&\|u\|^2_{\widetilde{L}^{\infty}([0,t];\,B^{s-1}_{2,2})}+c2\nu \|u\|^2_{\widetilde{L}^2([0,t];\,B^s_{2,2})}\\
\leq &(1+c2\nu t)M_0+C\int^t_0\|u\|_{L^{\infty}}\|u\|_{B^{s}_{2,2}}\|u\|_{
B^{s-1}_{2,2}}dt+C
\int^t_0\|b\|^2_{B^{s}_{2,2}}\|u\|_{B^{s-1}_{2,2}}dt'.
\end{align*}
Since $s>\frac d2$, we can choose a $r>2$ such that $s-1+\frac 2r\in (\frac d2,s)$. Then using the embedding:
\begin{equation}\label{embd}
\|f\|_{L^{\infty}}\leq C\|f\|_{B^{\theta}_{2,2}},\quad \theta>d/2,
\end{equation}
Young's inequality and the interpolation inequality in Lemma \ref{Besovpro} deduces that
\begin{align*}
&\|u\|^2_{\widetilde{L}^{\infty}([0,t];\,{B}^{s-1}_{2,2})}+c2\nu \|u\|^2_{\widetilde{L}^2([0,t];\,{B}^s_{2,2})}\\
\leq &2(1+2c\nu t)M_0+C\int^t_0\|u\|_{B^{s-1+\frac2r}_{2,2}}
\|u\|_{{B}^{s}_{2,2}}\|u\|_{B^{s-1}_{2,2}}dt^{\prime}+C
\int^t_0\|b\|^2_{B^{s}_{2,2}}\|u\|_{B^{s-1}_{2,2}}dt^{\prime}\\
\leq &2(1+2c\nu t)M_0+C\int^t_0\|u\|^{2-\frac2r}_{B^{s-1}_{2,2}}\|u\|^{1+\frac 2r}_{B^{s}_{2,2}}dt^{\prime}+C\int^t_0\|b\|^2_{B^{s}_{2,2}}
\|u\|_{B^{s-1}_{2,2}}dt^{\prime}\\
\leq &2(1+2c\nu t)M_0+c\nu\|u\|^2_{\widetilde{L}^2([0,t];\,
B^{s}_{2,2})}dt^{\prime}+\frac t2\|u\|^2_{\widetilde{L}^{\infty}([0,t];\,B^{s-1}_{2,2})}\\
&+C(\nu)\int^t_0
\|u\|^{4+\frac4{r-2}}_{\widetilde{L}^{\infty}([0,t^{\prime}];\,B^{s-1}_{2,2})}dt^{\prime}
+C\int^t_0\|b\|^4_{\widetilde{L}^{\infty}([0,t^{\prime}];\,B^{s}_{2,2})}
dt^{\prime}.
\end{align*}
Hence, we get that for any $t\leq1$,
\begin{equation}\label{u-s-1}
\begin{split}
&\|u(t)\|^2_{\widetilde{L}^{\infty}([0,t];\,{B}^{s-1}_{2,2})}+ \|u\|^2_{\widetilde{L}^2([0,t];\,{B}^s_{2,2})}\\
\leq &C_1(\nu)M_0+C_1(\nu)\int^t_0
\|u\|^{4+\frac4{r-2}}_{\widetilde{L}^{\infty}([0,t^{\prime}];\,B^{s-1}_{2,2})}dt^{\prime}
+C_2\int^t_0\|b\|^4_{\widetilde{L}^{\infty}([0,t^{\prime}];\,B^{s}_{2,2})}
dt^{\prime}.
\end{split}
\end{equation}
On the other hand, by Lemma \ref{pro.heat} and Proposition \ref{pro.comE}, we obtain that
\begin{equation*}
\begin{split}
\|u\|_{\widetilde{L}^1([0,T];\,\dot{B}^{s+1}_{2,2})}\leq & C(\nu)\|e^{\nu t\Delta}u_0\|_{\widetilde{L}^1_T(\dot{B}^{s+1}_{2,2})}+C(\nu)\|u\cdot\nabla u- b\cdot\nabla b\|_{\widetilde{L}^1_T(\dot{B}^{s-1}_{2,2})}\\
\leq &C(\nu)\|e^{\nu t\Delta}u_0\|_{\widetilde{L}^1_T(\dot{B}^{s+1}_{2,2})}+C(\nu)\|u\cdot\nabla u-b\cdot\nabla b \|_{L^1_T({B}^{s-1}_{2,2})}\\
\leq & C(\nu)\|e^{\nu t\Delta}u_0\|_{\widetilde{L}^1_T(\dot{B}^{s+1}_{2,2})}+C(\nu)T^{\frac{r-2}{2r}}\|u\|_{L^r([0,T];L^{\infty})}\|u\|_{\widetilde{L}^2_T(
\dot{B}^{s}_{2,2})}\\
&+C(\nu)\int^T_0\|b\|^2_{\dot{B}^{s}_{2,2}}\,dt.
\end{split}
\end{equation*}
According to the embedding that can be easily proved
\begin{equation}\label{embd1}
\|h\|_{L^q([0,T];\,L^{\infty})}\leq C\|h\|_{\widetilde{L}^q([0,T];\,B^{\sigma}_{2,2})},\quad\forall q\in[1,\infty],\quad  \sigma>d/2,
\end{equation} and the interpolation inequality in Lemma \ref{Besovpro}, we get that
\begin{equation*}
\begin{split}
&\|u\|_{\tilde{L}^1_T(\dot{B}^{s+1}_{2,2})}\\
\leq & C(\nu)\Big(\|e^{\nu t\Delta}u_0\|_{\widetilde{L}^1_T(\dot{B}^{s+1}_{2,2})}+T^{\frac{r-2}{2r}}\|u\|_{\widetilde{L}^r_T(
B^{s-1+\frac2r}_{2,2})}\|u\|_{\widetilde{L}^2_T(
{B}^{s}_{2,2})}+\int^T_0\|b\|^2_{\dot{B}^{s}_{2,2}}\,dt\Big)\\
\leq & C_2(\nu)\Big(\|e^{\nu t\Delta}u_0\|_{\widetilde{L}^1_T(\dot{B}^{s+1}_{2,2})}+T^{\frac{r-2}{2r}}\|u\|^{\frac{r-2}{r}}_{\widetilde{L}^{\infty}
_T(B^{s-1}_{2,2})}\|u\|^{1+\frac{2}{r}}_{\widetilde{L}^2_T(
{B}^{s}_{2,2})}+\int^T_0\|b\|^2_{\dot{B}^{s}_{2,2}}\,dt\Big).
\end{split}
\end{equation*}
This, together with \eqref{L^2}, gives that
\begin{equation}\label{u-s+1}
\begin{split}
\|u\|_{\tilde{L}^1_T(B^{s+1}_{2,2})}\leq & C_2(\nu)\|e^{\nu t\Delta}u_0\|_{\widetilde{L}^1_T(\dot{B}^{s+1}_{2,2})}+TM_0^{\frac{1}2}
+C_2(\nu)\int^T_0\|b\|^2_{B^{s}_{2,2}}\,dt\\
&+C_2(\nu)T^{\frac{r-2}{2r}}\|u\|^{\frac{r-2}{r}}_{\widetilde{L}^{\infty}
_T(B^{s-1}_{2,2})}\|u\|^{1+\frac{2}{r}}_{\widetilde{L}^2_T(
B^{s}_{2,2})}.
\end{split}
\end{equation}
\indent Now, let's deal with the  term $\|b\|^2_{B^{s}_{2,2}}$. Taking the $L^2$ inner product of the second equation of \eqref{eq.loc} with $\Delta_j b$ on $\RR^d,$  by H\"{o}lder's inequality and the fact that
\[\int_{\RR^d}(u\cdot \nabla){\Delta}_j b\cdot{\Delta}_j b\,dx=0,\]
we have that
\begin{equation}\label{bL2}
\begin{split}
\frac12\frac{d}{dt}\|{\Delta}_j b\|^2_{L^2}
\leq &\|[{\Delta}_j,\,u\cdot \nabla]b\|_{L^2}\|{\Delta}_j b\|_{L^2}+\|{\Delta}_j(b\cdot \nabla u)\|_{L^2}\|{\Delta}_j b\|_{L^2}\|{\Delta}_j b\|_{L^2}.
\end{split}
\end{equation}
By Lemma \ref{Dancom01}, Proposition \ref{Besovpro} and the embedding relationship ${B}^{s}_{2,2}\hookrightarrow{B}^{\frac d2}_{2,1}$, we obtain that
\begin{equation}\label{u-b}
\begin{split}
&\int_0^t\|[{\Delta}_j,\,u\cdot \nabla]b\|_{L^2}+\|{\Delta}_j(b\cdot \nabla u)\|_{L^2}\,dt'\\
\leq &Cc_j2^{-js}\int_0^t\|\nabla u\|_{{B}^{\frac d2}_{2,\infty}\cap L^{\infty}}\|b\|_{{B}^{s}_{2,2}}\,dt\\
\leq&Cc_j2^{-js}\|\nabla u\|_{\widetilde L^1([0,t];{B}^{\frac d2}_{2,1})}\|b\|_{\widetilde L^{\infty}([0,t]{B}^{s}_{2,2})}\\
\leq &Cc_j2^{-js}\|u\|_{\widetilde L^1([0,t];{B}^{s+1}_{2,2})}\|b\|_{\widetilde L^{\infty}([0,t]{B}^{s}_{2,2})}.
\end{split}
\end{equation}
Therefore, dividing both sides of \eqref{bL2} by $\|\dot{\Delta}_j b\|_{L^2}$, taking $L^1([0,t])$-norm, plugging \eqref{u-b} into \eqref{bL2},  we obtain that
\begin{align*}
\|{\Delta}_j b(t)\|_{\widetilde{L}^{\infty}_t(\,L^2)}\leq \|{\Delta}_j b_0\|_{L^2}+Cc_j2^{-js}
\|u\|_{\widetilde{L}^{1}_t({B}^{s+1}_{2,2})}\|b\|_{\widetilde L^{\infty}_t({B}^{s}_{2,2})}.
\end{align*}
Multiplying both sides of the above inequality by $2^{js}$ and  taking $\ell^2(j\ge-1)$-norm deduces that
\begin{equation}\label{b-s}
\begin{split}
\|b\|_{\widetilde{L}^{\infty}_t({B}^s_{2,2})}\leq \|b_0\|_{{B}^s_{2,2}}+C_3\|u\|_{\widetilde{L}^{1}_t(B^{s+1}_{2,2})}
\|b\|_{\widetilde{L}^{\infty}_t(B^s_{2,2})}.
\end{split}
\end{equation}
\indent Next, we will show that there exists a $T_{\ast}$ such that
$\|b\|_{\widetilde{L}^{\infty}([0,t];\,{B}^{s}_{2,2})}\leq 2\|b_0\|_{{B}^{ s}_{2,2}}$ for any $t\in [0,T_{\ast}]$. Notice that $$\|e^{\nu t\Delta}u_0\|_{\widetilde{L}^1_{T}(\dot{B}^{s+1}_{2,2})}\leq \sum_j 2^{2j(s-1)}\|\dot\Delta_ju_0\|_{L^2}(1-e^{-C2^{2j}T})\leq\|u_0\|_{\dot{B}^{s-1}_{2,2}}, $$ set
\[T'\triangleq\sup\Big\{T\in [0,T_*]\,\Big|\,\|b\|_{\widetilde{L}^{\infty}_T(\dot{B}^{s}_{2,2})}\leq 2\|b_0\|_{\dot{B}^{ s}_{2,2}}\Big\}\]
where $T_{\ast}$ satisfies
\begin{equation}\label{T1}
\Big(1-\frac{r}{r-2}C_1(\nu)T_{\ast}\big(C_{1}(\nu)M_0+16C_2\|b_0\|^4_{B^s_{2,2}}T_{\ast}\big)
\Big)^{-\frac{r-2}{r}}<2.
\end{equation}
and
\begin{equation}\label{T2}
\begin{split}
&C_2(\nu)\|e^{\nu t\Delta}u_0\|_{\widetilde{L}^1_{T'}(\dot{B}^{s+1}_{2,2})}+T'\big(M_0^{\frac{1}2}
+C_2(\nu)M_0\big)\\
&\quad+C_2(\nu)(T')^{\frac{r-2}{2r}}(2C_1(\nu)M_0+32C_2M^4_0)<\frac{1}{2C_3}.
\end{split}
\end{equation}
Suppose that $T'<T_{\ast}$,  we get that from \eqref{u-s-1}
\begin{align*}
&\|u\|^2_{\widetilde{L}^{\infty}([0,T'];\,{B}^{s-1}_{2,2})}+c2\nu \|u\|^2_{\widetilde{L}^2([0,T'];\,{B}^s_{2,2})}\\
\leq &C_1(\nu)M_0+C_1(\nu)\int^{T'}_0
\|u\|^{4+\frac4{r-2}}_{\widetilde{L}^{\infty}([0,t];\,B^{s-1}_{2,2})}dt
+16C_2T'\|b_0\|^4_{B^{s}_{2,2}}.
\end{align*}
Applying Lemma \ref{lem.gron} to the above inequality, from \eqref{T1}, we have that
\begin{equation}\label{u-L2'}
\|u\|^2_{\widetilde{L}^{\infty}_{T'}({B}^{s-1}_{2,2})}+c2\nu \|u\|^2_{\widetilde{L}^2_{T'}({B}^s_{2,2})}
\leq 2C_1(\nu)M_0++32C_2T\|b_0\|^4_{B^{s}_{2,2}}
\end{equation}
Substituting \eqref{u-L2'} into \eqref{u-s+1}, we get that
\begin{align*}
\|u\|_{\widetilde{L}^1_{T'}(B^{s+1}_{2,2})}\leq &C_2(\nu)\|e^{\nu t\Delta}u_0\|_{\widetilde{L}^1_{T'}(\dot{B}^{s+1}_{2,2})}+T'M_0^{\frac{1}2}
+C_2(\nu)T'\|b_0\|^2_{B^{s}_{2,2}}\\
&+C_2(\nu)(T')^{\frac{r-2}{2r}}(2C_1(\nu)M_0+32C_2T'\|b_0\|^4_{B^{s}_{2,2}})\\
\leq &C_2(\nu)\|e^{\nu t\Delta}u_0\|_{\widetilde{L}^1_{T'}(\dot{B}^{s+1}_{2,2})}+T'M_0^{\frac{1}2}
+C_2(\nu)T'M_0\\
&+C_2(\nu)(T')^{\frac{r-2}{2r}}(2C_1(\nu)M_0+32C_2M^4_0)<\frac{1}{2C_3}.
\end{align*}
This, together with \eqref{L^2}, implies that
\[\|b\|_{\widetilde{L}^{\infty}_T(B^{s}_{2,2})}<2\|b_0\|_{B^s_{2,2}},\]
contradicting the maximality of $T$. Hence $T=T_*$.
This fact, together with \eqref{u-s-1}, \eqref{u-s+1} and \eqref{b-s}, entails that the require results \ref{pro-pri}.
\end{proof}
\subsection{Convergence of the solution sequences}\label{CON}
\begin{proposition}\label{converge}
The solutions $(u_n,b_n)$ of approximate system \eqref{eq.NMHDAPp} is Cauchy with respect to $n$  in $C([0,T_{\ast});L^2(\RR^d))\times C([0,T_{\ast});L^2(\RR^d)).$ Moreover, the limits satisfy with that $$u\in C([0,T_{\ast});B_{2,2}^{s-1}(\RR^d))\cap L^2([0,T_{\ast});B_{2,2}^{s}(\RR^d))\cap L^1([0,T_{\ast});B_{2,2}^{s+1}(\RR^d)),$$ $$b\in C([0,T_{\ast});B_{2,2}^{s}(\RR^d)).$$
\end{proposition}
\begin{proof}
  We firstly prove the solutions $(u_n,b_n)$ of the approximate system $\eqref{eq.NMHDAPp}$ is Cauchy with respect to $n$. Assume $(u_n,b_n),(u_p,_p)$
are  any two  solutions of the approximate system $\eqref{eq.NMHDAPp}$.  They all satisfy with Proposition \ref{bound}:
$$\|u_n\|^2_{C([0,T_{\ast}];B^{s-1}_{2,2})}+\|u_n\|_{\widetilde L^1_{T_{\ast}} (B^{s+1}_{2,2})}+\|u_n\|^2_{\widetilde L^2_{T_{\ast}} (B^{s}_{2,2})}+\|b_{n}\|^2_{C([0,T_{\ast}];B^{s}_{2,2}}\le C_{\ast}.$$

The differences $u_n-u_p,\;b_{n}-b_{p}$ (suppose $p>n$) satisfy the following system:
\begin{align*}
\begin{cases}
\partial_t(u_n-u_p)-\nu\Delta(u_n-u_p)=J_n\PP(b_n\cdot\nabla b_n)-J_p\PP(b_p\cdot\nabla b_p)\\
\quad\quad\quad\quad\quad\quad-J_n\PP[(u_n\cdot\nabla)u_n]+J_p\PP[(u_p\cdot\nabla)u_p],\\
\partial_t(b_{n}-b_{p})=-J_{n}(u_n\cdot\nabla b_n)+J_{p}(b_p\cdot\nabla u_p)
\end{cases}
\end{align*}
\indent Applying localization operator $\dot\Delta_j$ to both sides of the above system, by the $L^2$ energy estimate, divergence free condition and integrating on $[0,t]$, we obtain that
\begin{align*}
 &\|\dot\Delta_j(u_n-u_p)\|^2_{L^2}+\|\dot\Delta_j(b_{n}-b_{p})\|^2_{L^2}+ \nu\int_0^t\|\dot\Delta_j\nabla (u_n-u_p)\|_{L^2}^2\\
=&\|\dot\Delta_j(u_{0n}-u_{0p})\|^2_{L^2}+\|\dot\Delta_j(b_{0n}-b_{0p})\|^2_{L^2}\\&+\int_0^t\langle \dot\Delta_j(J_n\PP(b_n\cdot\nabla b_n)
-J_p\PP(b_p\cdot\nabla b_p)),\dot\Delta_j(u_n-u_p)\rangle\\&-\langle (J_n\PP[(u_n\cdot\nabla,\dot\Delta_j]u_n-J_p\PP[(u_p\cdot\nabla,\dot\Delta_j] u_p,\dot\Delta_j(u_n-u_p)\rangle\\
&+J_{n}\dot\Delta_j((b_n\cdot\nabla)u_{n})-J_{p}\dot\Delta_j((b_p\cdot\nabla)u_p),\dot\Delta_j(b_{n}-b_{p})\rangle\\
&-J_{n}\dot\Delta_j((u_n\cdot\nabla)b_{n})-J_{p}\dot\Delta_j((u_p\cdot\nabla)b_p),\dot\Delta_j(b_{n}-b_{p})\rangle\,dt'
=:\sum\nolimits_{i=0}^4 E_i.
\end{align*}
where $E_0=\|\dot\Delta_j(u_{0n}-u_{0p})\|^2_{L^2}+\|\dot\Delta_j(b_{0n}-b_{0p})\|^2_{L^2}$.\\
\indent Split each $E_i\;(i=1,2,3)$ into three parts. We only deal with the two more difficult terms: $E_3$ and  $E_4$.
\begin{align*}
E_3=&\int_0^t\langle J_{n}\dot\Delta_j(b_n\cdot\nabla u_{n})-J_{p}\dot\Delta_j(b_p\cdot\nabla u_{p}),\dot\Delta_j(b_{n}-b_{p})\rangle dt'\\
=&\int_0^t\langle (J_{n}-J_{p})\dot\Delta_j(b_n\cdot\nabla u_{n}),\dot\Delta_j(b_{n}-b_{p})\rangle\\
&-\langle J_{p}\dot\Delta_j((b_n-b_p)\cdot\nabla u_{n}),\dot\Delta_j(b_{n}-b_{p})\rangle\\
&-\langle J_{p}\dot\Delta_j(b_p\cdot \nabla (u_{n}-u_p)),\dot\Delta_j(b_{n}-b_{p})\rangle dt'=:\sum\nolimits_{i=1}^3 E_{3i}.
\end{align*}
For $E_{31}$,  summing up over $j\in\ZZ$,  we obtain that by H\"older's inequality
\begin{equation}\label{E-31}
\begin{split}
&\sum_j|E_{31}|\leq\sum_j\int_0^t\|(J_{n}-J_{p})\dot\Delta_j(b_n\cdot\nabla u_{n})\|_{L^2}\|\dot\Delta_j(b_{n}-b_{p})\|_{L^2}\\
\leq &\Big\{\int_0^t\|(J_{n}-J_{p})\dot\Delta_j(b_n\cdot\nabla u_{n})\|_{L^2}dt'\Big\}_{\ell^2(\ZZ)}\|b_{n}-b_{p}\|_{\widetilde L_t^{\infty}(\dot B^0_{2,2})}\\
\leq &\frac{C}{n^{\epsilon}}\Big\{2^{j\epsilon}\int_0^t\|\dot\Delta_j(b_{n}\cdot\nabla u_n)\|_{L^2}dt'\Big\}_{\ell^2(\ZZ)}\|b_{n}-b_{p}\|_{\widetilde L_t^{\infty}(\dot B^0_{2,2})}\\
\le&\frac{C}{n^{\epsilon}}\|b_{n}\cdot\nabla u_n\|_{\widetilde L^1_t(\dot B^{\epsilon}_{2,2})}\|b_{n}-b_{p}\|_{\widetilde L_t^{\infty}(\dot B^0_{2,2})}\\
\leq& \frac{C}{n^{\epsilon}}\|u_n\|_{\widetilde L^{1}_t(B^{s+1}_{2,2})}\|b_{n}\|_{\widetilde L^{\infty}_t(B^{s}_{2,2})}\|b_{n}-b_{p}\|_{\widetilde L_t^{\infty}(\dot B^0_{2,2})}
\end{split}
\end{equation}
where  we choose $\epsilon$ satisfying $0<\epsilon<s-\frac d2$ and  use the embedding that $B^{s}_{2,2}\hookrightarrow B^{\epsilon}_{2,2}\hookrightarrow \dot B^{\epsilon}_{2,2}$ together with Lemma \ref{Besovpro}.\\
$E_{32}$ and $E_{33}$ can be estimated similarly as follows
\begin{equation*}
\begin{split}
\sum_j\sum ^3_{k=2}|E_{3k}|\leq &\frac 12\|\nabla u_n-\nabla u_p\|^2_{\widetilde L^{2}_t(\dot B^{0}_{2,2})}+C\Big(\| u_p\|_{\widetilde L^{1}_t(B^{s+1}_{2,2})}+t\|b_{n}\|^2_{\widetilde L^{\infty}_t(B^{s}_{2,2})}\Big)\|b_{n}-b_{p}\|^2_{\widetilde L^{\infty}_t(\dot B^{0}_{2,2})}.
\end{split}
\end{equation*}
where we use the Young's inequality.
Then, we have  that
\begin{equation*}
\begin{split}
&\sum_j\sum ^3_{k=1}|E_{3k}|
\leq C\int_0^t\|\nabla u_n-\nabla u_p\|_{\dot B^{0}_{2,2}}\|b_{n}\|_{B^{s}_{2,2}}\|b_{n}-b_{p}\|_{\widetilde L^{\infty}([0,t'];\dot B^{0}_{2,2})}dt'
\\&\quad\quad\quad\quad+(\|\nabla u_p\|_{\widetilde L^{1}([0,t'];B^{s}_{2,2})}
+t\|b_{n}\|_{\widetilde L^{\infty}([0,t'];B^{s}_{2,2})})\|b_{n}-b_{p}\|^2_{\widetilde L^{\infty}([0,t'];\dot B^{0}_{2,2})}\\
&\le\frac{C}{n^{\epsilon}}(\|u_n\|_{\widetilde L^{1}([0,t];B^{s+1}_{2,2})}+t\|b_{n}\|_{\widetilde L^{\infty}([0,t];B^{s}_{2,2})})\|b_{n}\|_{\widetilde L^{\infty}([0,t];B^{s}_{2,2})})+\frac 12\|\nabla u_n-\nabla u_p\|^2_{\widetilde L^{2}([0,t'];\dot B^{0}_{2,2})}\\&\quad\quad\quad\quad+C(\| u_p\|_{\widetilde L^{1}([0,t'];B^{s+1}_{2,2})}
+t\|b_{n}\|^2_{\widetilde L^{\infty}([0,t'];B^{s}_{2,2})})\|b_{n}-b_{p}\|^2_{\widetilde L^{\infty}([0,t'];\dot B^{0}_{2,2})}.
\end{split}
\end{equation*}
For $E_4$, we decompose into
\begin{align*}
E_4=&-\int_0^t\langle J_{n}\dot\Delta_j(u_n\cdot\nabla b_{n})-J_{p}(u_p\cdot\nabla b_{p}),\dot\Delta_j(b_{n}-b_{p})\rangle dt'\\
=&-\int_0^t\langle (J_{n}-J_{p})\dot\Delta_j(u_n\cdot\nabla b_{n}),\dot\Delta_j(b_{n}-b_{p})\rangle\\
&-\langle J_{p}\dot\Delta_j[(u_n-u_p)\cdot\nabla b_{n}],\dot\Delta_j(b_{n}-b_{p})\rangle\\
&-\langle J_{p}\dot\Delta_j[u_p\cdot\nabla(b_{n}-b_{p})],\dot\Delta_j(b_{n}-b_{p})\rangle dt'
=:\sum\nolimits_{i=1}^3 E_{4i}
\end{align*}
Just as $E_{31}$, we have
\begin{eqnarray*}
\sum_j|E_{41}|\le \frac{1}{n^{\epsilon}}\| u_n\|_{\widetilde L^{\infty}_t(B^{s+1}_{2,2}}\|b_{n}\|_{\widetilde L^{\infty}_t(B^{s}_{2,2}}\|b_{n}-b_{p}\|_{\widetilde L^{\infty}_t(\dot B^{0}_{2,2})}.
\end{eqnarray*}
Before dealing with $E_{42}$, we need the following estimates:\\
when $d=2$, $0<\varepsilon<s-\frac{d}{2}$, we have$$\|fg\|_{L^2}\le C\|f\|_{L^{\frac 2{\varepsilon}}}\|g\|_{L^{\frac 2{1-\varepsilon}}}\le C\|f\|_{B^{1-\varepsilon}_{2,2}}\|g\|_{B^{\varepsilon}_{2,2}}\le C\|f\|_{B^1_{2,2}}\|g\|_{B^{s-1}_{2,2}};$$
when $d=3$, we have$$\|fg\|_{L^2}\le C\|fg\|_{B^0_{2,2}}\le C\|fg\|_{B^{\delta}_{2,1}}\leq C \|f\|_{B^1_{2,2}}\|g\|_{B^{\frac{d}{2}-1+\delta}_{2,2}}\leq C\|f\|_{B^1_{2,2}}\|g\|_{B^{s-1}_{2,2}}$$
by Lemma \ref{Besovpro} as $0<\delta<s-\frac{d}{2}$.
This implies
\begin{align*}
\sum_j|E_{4}|&\le C\int_0^t\|u_n-u_p\|_{B^1_{2,2}}\|\nabla b_{n}\|_{B^{s-1}_{2,2}}\|b_{n}-b_{p}\|_{\widetilde L^{\infty}_t(\dot B^{0}_{2,2})}dt'\\
&\le \frac{1}{4}\|u_n-u_p\|^2_{\widetilde L^{2}_t(\dot B^{0}_{2,2})}+Ct\|b_{n}\|^2_{\widetilde L^{\infty}_t(B^s_{2,2})}\|b_{n}-_{p}\|^2_{\widetilde L^{\infty}_t(\dot B^{0}_{2,2})}.
\end{align*}
Noting that fact  $J_pJ_n=J_n$ when $p>n$,  we get $E_{43}=0$  by integrating by parts and using the divergence free condition.\\ In addition,  $E_1,\;E_2$ can be easily estimated similarly.\\
Plugging the above estimates together with estimates for $E_1,\;E_2$ and using Young's inequality, we obtain that
\begin{align*}
&\|u_n-u_p\|^2_{\widetilde L^{\infty}_t(\dot B^{0}_{2,2})}
+\|b_{n}-b_{p}\|^2_{\widetilde L^{\infty}_t(\dot B^{0}_{2,2})}+ \|\nabla (u_n-u_p)\|^2_{\widetilde L^{2}_t(\dot B^{0}_{2,2})}\\
\le&\frac{C}{n^{\epsilon}}\Big(((\|u_0\|^2_{B^{s-1}_{2,2}}+\|b_0\|^2_{B^{s}_{2,2}})
+t+(\|u_n\|^2_{\widetilde L^{\infty}_tB^{s}_{2,2}}+\| u_n\|_{\widetilde L^{1}_tB^{s+1}_{2,2}}+\|b_{n}\|^2_{\widetilde L^{\infty}_{t}(B^{s}_{2,2})})\\ &\times
(\|u_n-u_p\|_{\widetilde L^{\infty}_{t}(\dot B^{0}_{2,2})}+\|b_{n}-b_{p}\|_{\widetilde L^{\infty}_{t}(\dot B^{0}_{2,2})})\Big)\\
&+ \widetilde C(t+\| u_p\|_{\widetilde L^{1}_tB^{s}_{2,2}}+t\|b_{n}\|^2_{\widetilde L^{\infty}_tB^{s}_{2,2}})
(\|u_n-u_p\|^2_{\widetilde L^{\infty}_{t}(\dot B^{0}_{2,2})}+\|b_{n}-b_{p}\|^2_{\widetilde L^{\infty}_{t}(\dot B^{0}_{2,2})}).
\end{align*}
Denote
$$Y(t)=\|u_n-u_p\|^2_{\widetilde L^{\infty}_t(\dot B^{0}_{2,2})}+\|b_{n}-b_{p}\|^2_{\widetilde L^{\infty}_t(\dot B^{0}_{2,2})}.$$
By Proposition \ref{bound} and let $t$ small enough such that $\widetilde C(t+\| u_p\|_{\widetilde L^{1}_tB^{s}_{2,2}}+t\|b_{n}\|^2_{\widetilde L^{\infty}_tB^{s}_{2,2}})
\leq \frac 12$, we get that
$$Y(t)\le \frac{CC_{\ast}}{n^{\epsilon}}\longrightarrow0, n\longrightarrow+{\infty}.
$$
 Hence, $\{(u_n,b_n)\}_{n}$ is Cauchy for $n$ in $C([0,T_{\ast});L^2(\RR^d)\times L^2(\RR^d))$, by interpolation, in $C([0,T_{\ast});B^{s'-1}_{2,2})\times C([0,T_{\ast});B^{s'}_{2,2})$ for any $s'<s$. The limit $(u,b)$ is in $\widetilde L^{\infty}([0,T_{\ast});B^{s-1}_{2,2})\times \widetilde L^{\infty}([0,T_{\ast});B^{s}_{2,2})$ by using Fatou property for the Besov spaces. This, together with Proposition \ref{bound}
passing to the limit, deduces that by Fatou lemma $$u\in C([0,T_{\ast});B_{2,2}^{s-1}(\RR^d))\cap L^2([0,T_{\ast});B_{2,2}^{s}(\RR^d))\cap \widetilde L^1([0,T_{\ast});B_{2,2}^{s+1}(\RR^d)),$$ and $$b\in C([0,T_{\ast})B_{2,2}^{s}(\RR^d)).$$ Then we have proved Proposition \ref{converge}.
\end{proof}

\subsection{Uniqueness}\label{Uniquensee}
\begin{proposition}\label{Uniqueness}
The solution $(u,b)$ of system \eqref{eq.NMHD} in the previous step is unique.
\end{proposition}
\begin{proof}
Assume $(u_1,b_1),(u_2,b_2)$
are two any solutions of the system $\eqref{eq.NMHD}$. Certainly, they all satisfy with Proposition \ref{pro-pri}. Denote $u=u_1-u_2,\,b=b_1-b_2,\,\pi=\pi_1-\pi_2$, then the difference satisfies the following system:
\begin{equation}\label{U0}
\left\{\begin{array}{ll}
\partial_t u-\nu\Delta u+u_1\cdot\nabla u+u\cdot\nabla u_2+\nabla \pi=b_1\cdot\nabla b+b\cdot\nabla b_2,\\
\partial_t b+u_1\cdot\nabla b+u\cdot\nabla b_2=b_1\cdot\nabla u+b\cdot\nabla u_2,\\
\Div u=\Div b=0,\\
(u,b)|_{t=0}=(0,0).
\end{array}\right.
\end{equation}
The proof sketch is very similar to prove Proposition \ref{converge}, we omit it and  finally have small $t_1>0$ enough such that for any $0<t\leq t_1$
$$\|u\|^2_{\widetilde L^{\infty}_t(\dot B^{0}_{2,2})}+\|b\|^2_{\widetilde L^{\infty}_t(\dot B^{0}_{2,2})}\leq \frac 12 (\|u\|^2_{\widetilde L^{\infty}_t(\dot B^{0}_{2,2})}+\|b\|^2_{\widetilde L^{\infty}_t(\dot B^{0}_{2,2})}). $$
Therefore, we have $(u,b)=0$ on $[0,t_1]$. Using a continuity argument ensures that $(u_1,b_1)=(u_2,b_2)$ on $[0,T]$. This concludes the proof.
\end{proof}
Combining  Proposition \ref{pro-pri}, Proposition \ref{converge} and Proposition \ref{Uniqueness} together, we complete the proof of Theorem \ref{INTR}.

\end{document}